\magnification=\magstep1
\font\mittel=cmbx10 scaled \magstep1
\font\gross=cmbx10 scaled \magstep2
\font\dtsch=eufm10
\def\ni{\noindent}
\def\x{\times}
\def\.{\cdot}

\def\d{\partial}
\def\t{\vartheta}
\def\r{{\cal R}}
\def\g{{\hbox{{\dtsch g}}}}
\def\h{{\hbox{{\dtsch h}}}}
\def\k{{\hbox{{\dtsch k}}}}
\def\so{{\hbox{{\dtsch so}}}}
\def\u{{\hbox{{\dtsch u}}}}
\def\T{{\hbox{{\dtsch t}}}}
\def\p{{\hbox{{\dtsch p}}}}

\def\H{{\cal H}}
\def\V{{\cal V}}
\def\=*{\buildrel *\over =}

\def\Stroke#1{{\vrule height4.4pt width0.2pt} \kern-0.5pt{\scriptstyle{\rm #1}}}
\def\stroke#1{{\vrule height6.4pt width0.33pt} \kern-0.7pt{\rm #1}}
\def\mcr{\Stroke{R}}
\def\mcc{{\setbox0=\hbox{$\scriptstyle{\rm C}$}\raise0.056\ht0\hbox to0pt
{\hskip0.28\wd0\vrule height0.76\ht0 width0.24pt\hss}\copy0}}
\def\mcn{\Stroke{N}}
\def\cR{\stroke{R}}
\def\cC{{\setbox0=\hbox{\rm C}\raise0.055\ht0\hbox to0pt
{\hskip0.3\wd0\vrule height0.85\ht0 width0.35pt\hss}\copy0}}
\def\cN{{\vrule height6.4pt width0.3pt} \kern-0.7pt{\rm N}}
\def\R{\mathchoice{\cR}{\cR}{\mcr}{\mcr}}
\def\C{\mathchoice{\cC}{\cC}{\mcc}{\mcc}}
\def\N{\mathchoice{\cN}{\cN}{\mcn}{\mcn}}
\def\theorem#1{\par\medbreak\noindent\bf#1\quad\enspace\sl\nobreak}
\def\endtheorem{\par\medskip\rm}
\def\proof{\par\medbreak\noindent{\bf Proof}\quad\enspace}
\def\remark{\par\medbreak\noindent{\bf Remark}\quad\enspace}
\def\qed{\vbox{\hrule
  \hbox{\vrule\hbox to 5pt{\vbox to
5pt{\vfil}\hfil}\vrule}\hrule}}
\def\lit{\par\noindent
               \hangindent=3\parindent\ltextindent}
\def\ltextindent#1{\hbox to \hangindent{#1\hss}\ignorespaces}

\def\endproof{\unskip \nobreak \hskip0pt plus 1fill
\qquad \qed \par}

\def\3{ss}
\def\o{\circ}
\def\X{{\bar X}}

\def\Z{{\bar Z}}
\def\ox{\otimes}
\def\<{\langle}
\def\>{\rangle}
\def\^{\wedge}
\def\Hom{\mathop{\rm Hom}\nolimits}

\def\Span{\mathop{\rm Span\,}\nolimits}
\def\End{\mathop{\rm End}\nolimits}

\def\trace{\mathop{\rm trace\,}\nolimits}
\def\Aut{\mathop{\rm Aut}\nolimits}

\overfullrule=0pt

\centerline{\gross 
                    K\"ahler submanifolds}
\centerline{\gross
              with parallel pluri-mean curvature}

\bigskip
\centerline {F.E. Burstall\footnote{$^{1)}$}
{Partially supported by The European Contract Human Potential Programme,
Research Training Network HPRN-CT-2000-00101}, 
J.-H. Eschenburg\footnote{$^{2)}$}
    {partially supported by Volkswagenstiftung and DLR-CNPq}, 
M.J. Ferreira\footnote{$^{3)}$}
{Work supported by FCT, Praxis XXI, Feder, project Praxis/2/2.1/MAT/125/94}, 
R. Tribuzy$^{2)}$}

\bigskip
\bigskip
\par{\narrower\ni
{\bf Abstract:} {\it We investigate the local geometry of a class of
K\"ahler submanifolds $M \subset \R^n$ which generalize surfaces of
constant mean curvature. The role of the mean curvature vector is
played by the $(1,1)$-part (i.e. the $dz_id\bar z_j$-components) of
the second fundamental form $\alpha$, which we call the pluri-mean
curvature.  We show that these K\"ahler submanifolds are
characterized by the existence of an associated family of isometric
submanifolds with rotated second fundamental form. Of particular
interest is the isotropic case where this associated family is
trivial. We also investigate the properties of the corresponding
Gauss map which is pluriharmonic.}

\bigskip
\ni

\bigskip
\ni
\par}

\bigskip
\bigskip
\leftline {\mittel 1. Introduction}
\medskip

\noindent Some surfaces in 3-space admit isometric deformations which
change the shape of the surface while preserving the intrinsic
metric.  Even the principal curvatures may be preserved while the
principal curvature directions are rotated under the deformation;
this happens precisely if the surface has constant mean curvature
(``cmc"). The best known example is the deformation of the catenoid
into the helicoid which transforms the meridians and the equator of
the catenoid into the helicoid's ruling lines and axis, thus rotating
the principal curvature directions by $45^o$.

In the present paper we will to investigate submanifolds of 
higher dimension and codimension allowing similar deformations. 
The surface will be replaced by a simply connected $m$-dimensional 
complex manifold $M$ with an immersion $f : M \to \R^n$ such that 
the induced metric on $M$ is K\"ahlerian, i.e. the almost complex 
structure $J$ on $TM$ is orthogonal and parallel; we call these
{\it K\"ahler immersions} for short. Let $\alpha$ 
denote the second fundamental form of $f$ and rotate it by 
putting $\alpha_{\t}(x,y) = \alpha(\r_{\t}x,\r_{\t}y)$, where
$$
              \r_{\t} = \cos(\t)I + \sin(\t)J.
$$
When does there exist a family of isometric immersions $f_{\t} : 
M \to \R^n$ with second fundamental form $\alpha_{\t}$? We will 
see in Theorem 1 that this happens precisely if the bilinear form
$$
 \alpha^{(1,1)}(x,y) := {1\over2}(\alpha(x,y) + \alpha(Jx,Jy))
$$
is parallel with respect to the connections on the tangent and 
normal bundles. In the case of a surface ($m = 1$) we have 
$\alpha^{(1,1)}(x,y) = \<x,y\>\. \eta$, hence 
$\alpha^{(1,1)}$ is parallel if and only if the mean curvature vector 
$\eta = {1\over2}\trace\alpha$ is parallel also. 
This motivates us to call $\alpha^{(1,1)}$ the {\it 
pluri-mean curvature} of $f$; in fact, for any complex curve $C 
\subset M$ the restriction of $\alpha^{(1,1)}$ to $TC$ is again the 
metric multiplied by the mean curvature vector of the surface 
$f|_C$. But while surfaces with nonzero constant mean curvature 
can only have essential codimension 1 or 2 (cf. [Y]), there are 
interesting substantial examples in higher dimensions and
codimensions (cf. Ch. 7).

The main part of the paper is devoted to studying the relationship 
between a K\"ahler immersion $f : M \to \R^n$ with parallel 
pluri-mean curvature (``ppmc") and its Gauss map $\tau: M \to Gr$ 
where $\tau(p) = df_p(T_pM)$ and $Gr$ is the Grassmannian of 
$2m$-dimensional subspaces of $\R^n$. Just as in the case of cmc 
surfaces (cf. [RV]), ppmc submanifolds are characterized by the 
pluriharmonicity of their Gauss maps (Theorem 2). Pluriharmonic 
maps also admit an associated family of deformations, 
and in fact the deformed Gauss map is the Gauss map of the 
deformed immersion (Theorem 3).

The Gauss map $\tau$ of a K\"ahler immersion has a refinement $\tau'$
called the {\it complex Gauss map} which takes account of the complex
structure: for any $p \in M$ we put $\tau'(p) = df(T'_pM)$. (Here we
have extended $df_p$ complex linearly to $T^cM = TM \ox \C$ and used
the $J$-eigenspace decomposition $T^cM = T'M + T''M$ with $J = i$ on
$T'M$ and $J = -i$ on $T''M$.) The map $\tau'$ takes values in the
set $Z_1$ of {\it isotropic} complex $m$-dimensional subspaces $E
\subset \C^n$, i.e. the complex conjugate $\bar E$ is perpendicular
to $E$ with respect to the Hermitian inner product, or equivalently
$\<E,E\> = 0$ for the symmetric inner product $\<x,y\> = \sum x_jy_j$
on $\C^n$. This space $Z_1$ can be viewed as a flag manifold fibering
over $Gr$, and then $\tau'$ is a horizontal lift of $\tau$. We will
show that $\tau'$ is pluriharmonic if and only if $\tau$ is also and
hence if and only if $f$ is ppmc (Theorem 5). In fact we can
characterize the complex Gau\3 maps of ppmc immersions among the
pluriharmonic maps into $Z_1$ (Theorem 6).

Alternatively, $Z_1$ can also be viewed as a complex submanifold 
of the complex Grassmannian $Gc$ of $m$-planes in $\C^n$. We also 
study the composition $j \o \tau'$ for the inclusion $j : Z_1 \to 
Gc$. This map is pluriharmonic only for special ppmc immersions 
which we call {\it half isotropic} (Theorem 7). These contain two 
interesting subclasses, characterized also by properties of 
$\tau'$: the {\it pluri-minimal} ones with zero pluri-mean 
curvature  ($\tau'$ is holomorphic, Theorem 4) and the {\it 
isotropic} ones where the associated family is trivial 
($\tau$ and $j\o\tau'$ are isotropic, Theorems 9 and 10). 
The first of these results is well known for 
surfaces: a surface is minimal if and only if its (complex) Gauss 
map is holomorphic. The second result is not interesting for 
surfaces in 3-space: Isotropy would mean that each tangent 
vector is a principal curvature direction, hence the surface 
must be a round sphere or a plane. But there are interesting 
examples in higher dimension, among them the standard embeddings 
of Hermitian symmetric spaces (see Ch. 7). We need some facts
on flag manifolds which are known in principle but not explicitly
worked out; we shall prove these statements in an appendix.

Part of the work was done during visits of the second and the 
last named authors at ICTP, Trieste, and at the Mathematisches 
Forschungsinstitut Oberwolfach. They wish to express their thanks 
to these institutions for hospitality and to Volkswagenstiftung 
and DLR-CNPq for financial support.

\bigskip
\bigskip
\leftline{\mittel 2. Associated families of immersions}

\medskip
\noindent
Let $M$ be a K\"ahler manifold of complex dimension $m$; this is 
a $2m$-dimensional Riemannian manifold with a parallel and 
orthogonal almost complex structure $J$ on $TM$. Since our theory 
is entirely local, we do not need completeness of $M$, however at 
some points we will need simple connectivity. We consider an isometric
immersion $f : M \to \R^n$ (a {\it K\"ahler immersion}). 
Let $\alpha : TM \otimes TM 
\to N$ be the corresponding second fundamental form defined by 
$\alpha(X,Y) = (\d_X\d_Yf)^N$ where $N = Nf = df(TM)^\perp$ 
denotes the normal bundle of $f$. Consider the parallel rotations 
$\r_{\t} = \cos(\t)I + \sin(\t)J$ for any $\t \in \R$ and let 
$\alpha_{\t} : TM \otimes TM \to N$,
$$
       \alpha_{\t}(x,y) = \alpha(\r_{\t}x,\r_{\t}y). 
$$
An {\it associated family} for $f$ is roughly speaking a 
one-parameter family of isometric immersions $f_{\t} : M \to 
\R^n$ with second fundamental form $\alpha_{\t}$.\footnote{$^*)$}
    {This was called {\it weak associated family} in [ET2].} 
This is not quite correct since the second fundamental forms of 
the two immersions $f$ and $f_\t$ take values in different 
spaces, the normal bundles of $f$ and $f_\t$. More precisely, a 
one-parameter family $f_\t : M \to \R^n$ of isometric immersions  
will be called an {\it associated family} of $f$ if their second 
fundamental forms $\alpha_{f_\t}$ satisfy
$$
         \psi_\t(\alpha_{f_\t}(x,y)) = \alpha_\t(x,y)
                                      = \alpha(\r_\t x,\r_\t y) 
\leqno\ (1)
$$
for some parallel bundle isomorphism $\psi_\t : Nf_\t \to Nf$. 
Our first theorem below will show under which conditions such 
immersions exist.

We need some more notation. The complexified tangent bundle $T^cM 
= TM \otimes \C$ of a K\"ahler manifold $M$ splits as $T^cM = T' 
\oplus T''$ where the components are the parallel eigenbundles of 
the almost complex structure $J$ with $J = i$ on $T'$ and $J = -
i$ on $T''$. Vectors in $T'$ are also called (1,0)-vectors and 
those in $T'' = \overline{T'}$ are (0,1)-vectors. Let $\pi'(x) = 
{1\over2}(x-iJx)$ and $\pi''(x) = {1\over2}(x+iJx)$ be the 
projections onto these subbundles. Extending $\alpha$ complex 
linearly to the complexified tangent and normal bundles, we put
$$
\alpha^{(1,1)}(x,y) = \alpha(\pi'x,\pi''y) + \alpha(\pi''x,\pi'y)
              = {1\over2}(\alpha(x,y) + \alpha(Jx,Jy)).
\leqno\ (2)
$$
As explained in the introduction, $\alpha^{(1,1)}$ will be called the
{\it pluri-mean curvature}, and $f$ is called an immersion with
{\it parallel pluri-mean curvature} ({\it ppmc}) if this tensor 
is parallel with respect to the tangent and normal connections. 
The following theorem which was partially obtained in [FT] shows 
the relation to associated families.

\theorem{Theorem 1.} Let $f : M \to \R^n$ be a K\"ahler 
immersion. Then $f$ has an associated family if and only if it 
has parallel pluri-mean curvature.
\endtheorem

\proof 
We are using the existence theorem for submanifolds (cf [Sp]): 
Let $M$ be a $p$-dimensional Riemannian manifold and $N$ a 
$k$-dimensional vector bundle over $M$ with a connection $D^N$. 
Further let $\alpha \in \Hom(S^2TM,N)$ where $S^2TM$ denotes the 
symmetric tensor product of $TM$. Then there is an isometric 
immersion $f : M \to \R^{p+k}$ with normal bundle $N$ (up to a 
parallel vector bundle isometry) and second fundamental form 
$\alpha$ if and only if the submanifold equations of Gau\3, 
Codazzi and Ricci are satisfied.

Let us apply this to $\alpha_{\t}$. The {\it Gau\3 equation} is
$$
    \<R(x,y)v,w\> = \<\alpha_{\t}(x,w),\alpha_{\t}(y,v)\>
           - \<\alpha_{\t}(x,v),\alpha_{\t}(y,w)\>
\leqno\ (G_{\t})
$$
In fact, this equation follows from $(G_0)$, the Gau\3 equation 
of $f$. The easiest way to see this is to use the splitting $T^cM 
= T' + T''$. On $T'$ we have $\r_{\t} = e^{i\t}$ while $\r_{\t} = 
e^{-i\t}$ on $T''$. We may assume that $x,y,v,w \in T' \cup T''$. 
In all possible cases, the right hand side of $(G_\t)$ picks up a 
common factor $e^{ik\t}$ for some $k$. The left hand side is zero 
as soon $x,y$ or $v,w$ have the same type (both in $T'$ or both 
in $T''$). This holds on any K\"ahler manifold since $R(x,y)T' 
\subset T'$ and $\<T',T'\> = 0$, thus $\<R(x,y)T',T'\> = 0$ for 
all $x,y \in T^c$ (where we have extended the inner product 
complex linearly to $T^cM$). For these cases $(G_{\t})$ follows 
from $(G_0)$. In the remaining cases, two of the vectors 
$x,y,v,w$ are in $T'$ and the other two in $T''$, and thus 
$(G_{\t})$ is the same as $(G_0)$.

Next we consider the {\it Codazzi equation}:
$$
       (D_x\alpha_{\t})(y,z) = (D_y\alpha_{\t})(x,z).
\leqno\ (C_{\t})
$$
This follows from $(C_0)$ (the Codazzi equation of $f$) provided 
that $x,y$ have the same type. But if $x \in T'$ and $y \in T''$, 
we get different factors in front of the two sides of $(C_{\t})$. 
Thus $(C_{\t})$ follows from $(C_0)$ precisely if 
$(D_{T'}\alpha)(T'',T^c)$ vanishes, but by $(C_0)$, this is the 
same as $(D_{T^c}\alpha)(T',T'')$. Thus $(C_{\t})$ holds if and 
only if $\alpha^{(1,1)}$ is parallel. 

It remains to consider the {\it Ricci equation}. For any $\xi \in 
N$ let $A^{\t}_\xi$ be the symmetric endomorphism of $TM$ defined 
by
$$
       \<A^{\t}_\xi x,y\> = \<\alpha_{\t}(x,y),\xi\> = 
 \<\alpha(\r_{\t}x,\r_{\t}y),\xi\> = \<A_\xi \r_\t x,\r_\t y\>,
$$
hence $A^\t_\xi = \r_\t^{-1}A_\xi\r_\t$. Then the Ricci equation 
is
$$
      \<R^N(x,y)\xi,\eta\> = \<[A^\t_\xi,A^\t_\eta]x,y\>
                           = \<[A_\xi,A_\eta]\r_\t x,\r_\t y\>.
\leqno\ (R_\t)
$$
Again this equation follows from $(R_0)$, the Ricci equation for
$f$, provided that $x,y \in T' \cup T''$ are of different type. 
But if, say, both $x,y$ are in $T'$, the right hand side is
multiplied by $e^{2i\t}$. Hence $(R_\t)$ follows from $(R_0)$ if 
and only if $R^N(T',T') = 0$. (Note that the case $x,y \in T''$ 
arises just by complex conjugation.) But the subsequent Lemma 
shows that this is not a new condition; it follows also from 
$D\alpha^{(1,1)} = 0$. This finishes the proof of Theorem 1.
\endproof

\theorem{Lemma 1.} If a K\"ahler immersion $f : M \to \R^n$ has 
parallel pluri-mean curvature, then $R^N(T',T') = 0$.
\endtheorem

\proof
Let $N^o \subset N$ denote the image of $\alpha^{(1,1)}$; since 
$\alpha^{(1,1)}$ is parallel, $N^o$ is a parallel subbundle of 
$N$. Let $(N^o)^\perp \subset N$ be its orthogonal complement. For 
any $\xi \in (N^o)^\perp$ and $x \in T'$, $\bar y \in T''$ we have 
€$\<A_\xi x,\bar y\> = \<\alpha(x,\bar y),\xi\> \in 
\<N^o,(N^o)^\perp\> = 0$. Since $T'$ and $T''$ are isotropic subspaces, 
this implies $A_\xi(T') \subset T''$, and by complex conjugation 
we also get $A_\xi(T'') \subset T'$.

We have to show that
$
         \<R^N(x,y)\xi,\eta\> = \<[A_\xi,A_\eta]x,y\>.
$
vanishes for all $x,y \in T'$ and $\xi,\eta \in N$. It is 
sufficient to consider the following two cases:

\smallskip
\item{(a)} $\xi,\eta \in (N^o)^\perp$,
\item{(b)} $\xi \in N^o$ and $\eta \in N$ arbitrary. 

\smallskip\ni
In Case (a), both $A_\xi$ and $A_\eta$ interchange $T'$ 
and $T''$. Hence the commutator $[A_\xi,A_\eta]$ preserves $T'$ 
which by isotropy of $T'$ implies $\<[A_\xi,A_\eta]T',T'\> = 0$. 
Case (b) will follow from the following more general fact which 
is well known and easy to prove by twofold covariant 
differentiation:

\theorem{Sublemma.} Let $E,F$ be vector bundles with connections 
$D^E$ and $D^F$ over some smooth manifold $M$. Let $\beta : E \to 
F$ be a parallel homomorphism, i.e. $\beta(D^E_Xe) = 
D^F_X\beta(e)$ for any section $e$ of $E$. Then $R^F(x,y)\beta e 
= \beta(R^E(x,y)e)$.
\endproof
\endtheorem

\ni
We apply the Sublemma to $\beta := \alpha^{(1,1)} : T' \ox T'' \to 
N^o$. According to Case (b), we may assume $\xi = \alpha(u,\bar 
v)$ for some $u \in T'$ and $\bar v \in T''$. Since $N^o \subset 
N$ is parallel, we have
$$
        R^N(x,y)\xi = R^{N^o}(x,y)\beta(u\ox \bar v) = 
            \beta(R^{T'\ox T''}(x,y)(u\ox \bar v)) = 0,
$$
recall that $R^{T'\ox T''}(x,y)(u\ox \bar v) = (R(x,y)u)\ox \bar 
v + u \ox R(x,y)\bar v$, and $R(x,y) = 0$ for $x,y \in T'$ since 
$M$ is a K\"ahler manifold.
\endproof

\bigskip
\bigskip
\leftline{\mittel 3. The Gau\3 map}

\medskip
\noindent
Let $M$ be a $p$-dimensional smooth manifold, $f : M \to \R^n$ an 
immersion and $Gr$ the Grassmannian of $p$-dimensional linear 
subspaces in $\R^n$. The {\it Gau\3 map} $\tau : M \to Gr$ 
assigns to each $p \in M$ the subspace $\tau(p) = df_p(T_pM) 
€\subset \R^n$. We view $Gr$ as a submanifold in the vector space 
$S(n)$ of all symmetric real $n \x n$-matrices by replacing a 
linear subspace $E$ with the orthogonal projection onto $E$ (which 
we will also call $E$). Then the tangent space $T_EGr$ is the set 
$S(E,E^\perp)$ of all self adjoint linear maps on $\R^n$ sending 
$E$ to $E^\perp$ and vice versa; it can be naturally identified 
with $\Hom(E,E^\perp)$. A smooth map $\phi : M \to Gr$ can be 
viewed as a vector bundle over $M$. In fact $\phi$ and 
$\phi^\perp$ are subbundles of the trivial bundle $M \x \R^n$ and 
thus they inherit a natural connection which is differentiation 
on $\R^n$ followed by projection onto the fibre. We may view 
$\phi^*TGr = T_\phi Gr = \Hom(\phi, \phi^\perp)$, and the pull 
back connection on $\phi^*TGr$ is just the natural connection on 
$\Hom(\phi,\phi^\perp)$. The differential $d\phi : TM \to 
\phi^*TGr$ can be computed as follows: For any section $s$ of 
$\phi$ (i.e. $s : M \to \R^n$ with $s(p) \subset \phi(p)$ for all 
$p \in M$) we have
$$
(\d_x\phi)\.s = \d_x(\phi\.s) - \phi\.\d_xs = \phi^\perp\.\d_xs
\eqno\ (3)
$$
where $\phi(p)$ for any $p \in M$ is considered as a projection 
matrix on $\R^n$. In order to apply this to the Gau\3 map $\phi = 
\tau$, we use the section $s = df(Y) = \d_Yf$ where $Y$ is an 
arbitrary tangent vector field on $M$, and we obtain
$$
   (\d_X\tau)\.df(Y) = \tau^\perp(\d_X\d_Yf) = \alpha(X,Y).
\leqno\ (4)
$$

The following theorem due to [FT] generalizes the well known 
result of Ruh and Vilms [RV] which characterizes cmc surfaces by 
the harmonicity of their Gauss map. In higher dimension, 
harmonicity has to be replaced by pluriharmonicity: A smooth map 
$\phi : M \to S$ into a symmetric space $S$ is called {\it 
pluriharmonic} if its Levi form $Dd\phi^{(1,1)}$ (the restriction 
of the Hessian to $T' \ox T''$) vanishes. As always we view $d\phi$ as
a section of the bundle $\Hom(TM,f^*TS)$ with its natural connection
induced by the Levi-Civita connections on $M$ and $S$. In particular,
$d\tau$ is a section of $\Hom(TM,\Hom(\tau,\tau^\perp)) = \Hom(TM \otimes
\tau,N)$. Since $f : M \to \R^n$ is an isometric immersion, $df : TM \to
\tau$ is a parallel bundle isomorphism which will be used to identify
the bundles $TM$ and $\tau$. Using this identification and $(4)$ 
we have $d\tau = \alpha \in \Hom(TM \otimes TM,N)$ and $Dd\tau = D\alpha$.

\theorem{Theorem 2.} Let $M$ be a K\"ahler manifold and $f : M 
\to \R^n$ an isometric immersion. This has parallel pluri-mean 
curvature if and only if its Gau\3 map $\tau$ is pluriharmonic.
\endtheorem

\proof
$Dd\tau^{(1,1)} = 0$ if and only if 
for any $X \in T'$, $\bar Y \in T''$ and $W \in T^c$ we have $0 = 
(D_Xd\tau)(\bar Y)\.df(W) = (D_X\alpha)(\bar Y,W) = 
(D_W\alpha)(X,\bar Y)$, using Codazzi equation. Since $T'$ and $T''$
are parallel subbundles of $T^cM$, this is 
equivalent to $D(\alpha^{(1,1)}) = 0$.
\endproof

\medskip
The pluriharmonic map $\tau : M \to Gr$ has also an associated 
family: For any K\"ahler manifold $M$ and 
any symmetric space $S$, a family of smooth maps $\tau_\t : M \to S$ 
is called {\it associated} to $\tau = \tau_0$ if there is a 
€parallel bundle isomorphism $\phi_\t : \tau_\t^*TS \to \tau^*TS$ 
preserving the curvature tensor $R^S$ such that
$$
                \phi_\t \o d\tau_\t = d\tau \o \r_\t.
\leqno\ (5)
$$
It is known (cf. [ET2]) that a given smooth map $\tau : M \to S$ 
has a (unique) associated 
family if and only if it is pluriharmonic. 
We shall show next that the associated 
families of a ppmc immersion $f$ and its pluriharmonic Gau\3 map 
$\tau$ correspond to each other.

\theorem{Theorem 3.} Let $f : M \to \R^n$ be a ppmc immersion 
with Gau\3 map $\tau$ and let $f_\t$ be the associated family of 
$f$. Let $\tau_\t$ be the Gau\3 map of $f_\t$. Then $(\tau_\t)$ 
is the associated family of $\tau$.
\endtheorem

\proof
It suffices to show that the Gau\3 maps $\tau_\t$ of the 
immersions $f_\t$ form an associated family. Thus we have to find 
a parallel bundle map $\phi_\t : T_{\tau_\t}Gr \to T_\tau Gr$ 
satisfying $(5)$ above. Let $x,y \in T_pM$. On the one hand we 
have
$$
       d\tau(R_\t x) : df(y) \mapsto \alpha(\r_\t x,y),
$$
on the other hand
$$
          d\tau_\t(x) : df_\t(\r_{-\t}y) \mapsto 
             \alpha_\t(x,\r_{-\t}y) = \psi_\t(\alpha(\r_\t x,y)).
$$
Thus Equation $(5)$ is satisfied if for any $a \in 
T_{\tau_\t(p)}Gr = \Hom(\tau_\t(p),N_\t(p))$ we put
$$
            \phi_\t(a) = \psi_\t \o a \o \r_{-\t} \in   
               \Hom(\tau(p),N(p)) = T_{\tau(p)}Gr
$$
where we have identified both $\tau$ and $\tau_\t$ with $TM$ 
using $df$ and $df_\t$ and where $\psi_\t$ denotes the parallel 
isomorphism between the normal bundles (cf. $(1)$ in Ch.2). We 
see that $\phi_\t(p)$ acts by conjugating $a$ with the orthogonal 
$n\x n$-matrix $B$  mapping the subspaces $\tau_\t(p)$ and 
$N_\t(p)$ onto $\tau(p)$ and $N(p)$, with
$$
      B|_{\tau_\t(p)} = df_p \o R_\t \o (df_\t)_p^{-1},\ \ 
      B|_{N_\t(p)} = \psi_\t(p).
$$
Conjugation by $B \in O(n)$ is a global isometry on $Gr$ and thus 
preserves the curvature tensor of $Gr$. Moreover, $\phi_\t$ is 
parallel since so are $\psi_\t$ and $\r_\t$ as well as $df : TM 
\to \tau$ and $df_\t : TM \to \tau_\t$. Thus $\tau_\t$ is the 
associated family of $\tau$.
\endproof

\bigskip
\bigskip
\leftline{\mittel 4. The complex Gau\3 map}

\medskip
\noindent
The Gau\3 map $\tau$ of a K\"ahler manifold immersion $f : M \to 
\R^n$ records only the tangent planes without taking care of 
the complex structure. Therefore we introduce a refinement, the 
{\it complex Gau\3 map} $\tau'$. It takes values in the set $Z_1$ 
of all $m$-dimensional linear subspaces $E \subset \C^n$ which 
are {\it isotropic}, i.e. the bilinear inner product $\<x,y\> = 
\sum_j x_jy_j$ on $\C^n$ vanishes on $E \x E$. In fact we let 
$\tau' : M \to Z_1$,
$$
   \tau'(p) = df(T'_p) = \{df(x) - i\.df(Jx);\ x \in T_pM\}
                         \subset \C^n.
$$

The manifold $Z_1$ can be viewed in two different ways. On the 
one hand, it is a complex submanifold of the complex Grassmannian 
$Gc = G_m(\C^n)$ of all complex $m$-planes in $\C^n$. In fact, 
the complex structure on $Gc$ is induced by the complex Lie group 
$GL(n,\C)$ acting transitively on $Gc$, 
and $Z_1 \subset Gc$ is an orbit of the complex 
subgroup $O(n,\C)$ inducing a complex structure on $Z_1$. 
On the other hand $Z_1$ can be considered also as a flag manifold 
fibering over the {\sl real} Grassmannian $Gr$ (cf. Appendix): To 
any $E \in Z_1$ we may assign the orthogonal \footnote{$^*)$}
        {The terms ``orthogonal" or ``perpendicular" in a complex 
        vector space are always related to the {\sl Hermitian}
        inner product $(x,y) = \<x,\overline y\>$.} 
decomposition (``flag'') 
$\C^n = E + N + \bar E$ where $N = (E+\bar E)^\perp$, and the 
projection $\pi : Z_1 \to Gr$ is given by $\pi(E) = E + \bar E$ 
(we view the subspaces of $\R^n$ as complex subspaces of $\C^n$ 
which are invariant under complex conjugation). In terms of coset 
spaces we have $Z_1 = O_n/(U_m\x O_k)$ where $k = n-2m$, and 
$\pi : Z_1 \to Gr = O_n/(O_{2m}\x O_k)$ is the canonical 
projection. This is a Riemannian submersion (up to a scaling 
factor) for any $O_n$-invariant metric on $Z_1$ since the 
horizontal space (the reductive complement of $\so_{2m} \oplus 
\so_k$ in the Lie algebra $\so_n$) 
is irreducible with respect to the isotropy 
group $U_m \x O_k$ of $Z_1$. As a further consequence, the notions
``horizontal'' and ``super-horizontal'' agree for $Z_1$ (cf.
Appendix).

If we take the second view point considering $Z_1$ as a flag 
manifold over $Gr$, we have to replace $\tau'$ by
$$
                    \tau_1 = (\tau',N,\tau'')
$$
where $\tau'' = \overline{\tau'}$ and $N = (\tau'+\tau'')^\perp$; 
this is the complexified normal bundle of the immersion $f$. 
Clearly, $\pi \o \tau_1 = \tau$.

\theorem{Lemma 2.} Let $f : M \to \R^n$ be a K\"ahler immersion 
with second fundamental form $\alpha$ and complex Gau\3 map 
$\tau' : M \to Z_1 \subset Gc$. Then we have for any $v \in TM$ 
and $x' \in T'$ (whence $df(x') \in \tau'$)
$$
                d\tau'(v).df(x') = \alpha(v,x'). 
\leqno\ (6)
$$
Consequently $\tau_1 = (\tau',N,\tau'')$ is a (super-)horizontal 
lift of the real Gau\3 map $\tau$.

\endtheorem

\proof
We first view $Z_1 \subset Gc$. We may identify $TM$ with $\tau$ 
and $T'$ with $\tau'$ using $df$. Since $(T')^\perp = T'' + N$, 
we have (as for the real Grassmannian) $d\tau'(v).x' = 
(\d_vX')^{(T')^\perp} = (\d_vX')^{T''} + (\d_vX')^N$ where $X'$ 
is a (1,0) vector field extending $x'$. But $(\d_vX')^{T''} = 
(D_vX')^{T''} = 0$ because $T'$ is parallel with respect to the 
Levi-Civita connection $D$ of $M$. Moreover $(\d_vX')^N = 
\alpha(v,x')$ which shows (6).

Now consider $Z_1$ as a flag manifold over $Gr$. Then Equation 
(6) shows that $d\tau_1(v) = (d\tau'(v),dN(v),d\tau''(v))$ is a 
superhorizontal vector since it maps $\tau'$ into the next following
space $N$; in other words, $d\tau_1(v).\tau'$ has no component in 
$\tau''$ (cf. Equation (A5) in the appendix).
\endproof

\remark
The proof shows that the horizontality of $\tau_1$ is just 
another expression for the parallelity of the almost complex 
structure $J$ on $M$.

\bigskip
The first occasion where the complex Gau\3 map turned out to be 
useful was the characterization of pluriminimal submanifolds by 
holomorphicity of $\tau'$ (cf. [RT]). A similar statement for 
$\tau$ would not even make sense.

\theorem{Theorem 4.} An isometric K\"ahler immersion $f : M \to 
\R^n$ is pluriminimal (i.e. has zero pluri-mean curvature) if and 
only if $\tau_1 : M \to Z_1$ is holomorphic.
\endtheorem

\proof
The map $\tau_1 = (\tau',N,\tau'')$ is holomorphic if and only if 
$d\tau_1$ maps $T' = T'M$ into $T'Z$ or, more precisely (using 
Lemma 2), into ${\cal H}_1'$. In other words (cf. Appendix), 
$d\tau_1(v')$ for $v' \in T'$ is a linear map sending $\tau'$ 
into $N$ (which is always true by Lemma 2) and $N$ into $\tau''$. 
The latter property says that for any $w'' \in T''$ 
and $\xi \in N$
$$
                 0 = -\<d\tau_1(v').\xi,w''\>  
                   = \<\xi,d\tau_1(v').w''\> 
                    = \<\xi,\alpha(v',w'')\>
$$
which means that $\alpha^{(1,1)} = 0$.
\endproof

\theorem{Theorem 5.} An isometric K\"ahler immersion $f : M \to 
\R^n$ is ppmc if and only if $\tau_1 : M \to Z_1$ is a 
(super)horizontal pluriharmonic map.
\endtheorem

\proof
By Lemma 2 the complex Gau\3 map $\tau_1$ of any K\"ahler immersion 
$f$ takes values in the (super)horizontal bundle ${\cal H}_1$. 
Moreover $f$ is ppmc if and only if its real Gau\3 map $\tau$ is 
pluriharmonic (cf. Theorem 2). But $\tau_1$ is a horizontal 
lift of $\tau$ with respect to the Riemannian submersion $\pi : 
Z_1 \to Gr$. This implies that pluriharmonicity for $\tau$ 
and $\tau_1$ are equivalent. In fact, $\tau$ is pluriharmonic 
if and only if for any two commuting vector fields $V' \in T'$ and $W'' \in 
T''$ we have $D_{W''}d\tau(V') = 0$. Since $d\tau_1(V')$ is the 
horizontal lift of $d\tau(V')$, this is equivalent to 
$D_{W''}d\tau_1(V') = 0$, see the subsequent Lemma 3 for details. 
\endproof

\theorem{Lemma 3.} Let $Z,S$ be Riemannian manifolds and $\pi : Z \to S$
a Riemannian submersion. Let $M$ be any manifold and $\tau_1 : M \to Z$
be a horizontal map, i.e. $d\tau_1(TM) \subset \H$ where $\H \subset TZ$ is the
horizontal subbundle. Consider the O'Neill tensor $A : \H \otimes \H \to \V$
(where $\V = \H^\perp \subset TZ$ is the vertical bundle) given by
$$
                A(X,Y) = [X,Y]^{\V} = 2(D_XY)^{\V}
$$
for horizontal vector fields $X,Y$. Then $\tau_1^*A = 0$, i.e.
$(D_Wd\tau_1(V))^{\V} = 0$ for any two vector fields $V,W$ on $M$.
\endtheorem

\proof
Let $V,W$ be local vector fields on $M$ with $[V,W] = 0$. Locally 
we can write $d\tau_1(V) = \sum_i v_i(X_i\o\tau_1)$ and $d\tau_1(W) = 
\sum_j w_j (X_j\o\tau_1)$ where $v_i,w_j$ are functions on $M$ and 
$X_1,...,X_n$ form a basis of horizontal vector fields on $Z$. Then
$A(d\tau_1(V),d\tau_1(W)) = \sum_{ij}v_iw_j A(X_i,X_j) \o\tau_1
= \sum_{ij} v_iw_j(D_{X_i}X_j - D_{X_j}X_i)^{\V}\o\tau_1 
\=* (D_Vd\tau_1(W) - D_Wd\tau_1(V))^{\V} = 0$, due to the symmetry
of the hessian $Dd\tau_1$; at $\,*\,$ we have used the identity 
$D_V(X_j\o\tau_1) = D_{d\tau_1(V)}X_j = \sum_i v_i (D_{X_i}X_j)\o\tau_1$ which
is a defining property of the induced connection on vector field along
$\tau_1$ and which implies $\sum_{ij} v_iw_j (D_{X_i}X_j)\o\tau_1 =
D_Vd\tau_1(W)$.
\endproof

\medskip

Now we can characterize all ppmc immersions with values in the 
unit sphere $S^{n-1} \subset \R^n$ by their complex Gau\3 map. In 
principle we are able to decide whether or not a given horizontal 
pluriharmonic map $\tau_1 : M \to Z_1$ is the complex Gau\3 map 
of a ppmc K\"ahler immersion:

\theorem{Theorem 6.} Let $M$ be a K\"ahler manifold. A horizontal 
pluriharmonic map $\tau_1 = (\tau',N,\tau''): M \to Z_1$ is the 
complex Gau\3 map of a ppmc K\"ahler immersion $f : M \to S^{n-1} 
\subset \R^n$ if and only if there exists a real section 
$f$ of $N$ (a smooth map $f : M \to \R^n$ with $f(p) \in N_p$ 
for all $p \in M$) such that $df(T') = \tau'$.
\endtheorem

\proof
Clearly, if $f : M \to S^{n-1}$ is a K\"ahler immersion, the 
position vector $f$ is always normal and hence a section of the 
normal bundle $N$ with $df(T') = \tau'$. Further, if $f$ is ppmc 
then $\tau_1 = (\tau',N,\tau'')$ is horizontal pluriharmonic by 
the previous theorem. Conversely, suppose that such a map $\tau_1 = 
(\tau',N,\tau'')$ and a real section $f$ of $N$ with 
$df(T') = \tau'$ are given. 
Since the values of $df$ are perpendicular to $N$, hence 
to $f$, we have $\<f,f\> = const \neq 0$, and we may assume that 
$f$ takes values in $S^{n-1}$. In order to show that it is a ppmc 
immersion, by Theorem 2 we have to prove only that the metric 
induced by $f$ on $M$ is K\"ahler for the given complex 
structure. In general this is true (cf. [ET1]) if and only if 

\smallskip
\item\item{(a)}  $df(T')$ is isotropic and
\item\item{(b)}  $ddf^{(1,1)}$ takes values in the normal bundle 
                 of $f$.
 
\smallskip\noindent
(a) is true since $df(T') = \tau'$ is isotropic by definition of 
$Z_1$, and (b) holds since $\tau'$ differentiates into $N$ by 
horizontality of $\tau_1$. More precisely, let $V'$ and 
$W''$ be commuting (1,0) and (0,1) vector fields. Then 
$s := \d_{V'}f$ is a section of $\tau'$, and hence 
$(\d_{W''}s)^{(\tau')^\perp} = (\d_{W''}\tau).s \in N$ (cf. (3) 
in Ch. 3). Hence $\d_{W''}\d_{V'}f \in \tau' + N$. Similar we 
obtain $\d_{V'}\d_{W''}f \in \tau'' + N$. Since the two 
expressions agree, they must be contained in the intersection of 
the two bundles which is $N$.
\endproof

\medskip
Returning to the first view point $Z_1 \subset Gc$ we may ask if also
$\tau' : M \to Gc$ is pluriharmonic when $f$ is ppmc. In [FT] 
it was shown that an extra condition is needed: Let $N^o \subset 
N$ be the parallel subbundle spanned by the values of 
$\alpha^{(1,1)}$ and $N^1$ its orthogonal complement in $N$. The 
ppmc immersion $f$ is called {\it half isotropic} if 
$\alpha(T',T') \subset N^1$. The reason for this notation will 
become clear in the next chapter.

\theorem{Theorem 7.} Let  $M$ be a K\"ahler manifold and $f : M 
\to \R^n$ an isometric immersion with complex Gau\3 map $\tau' : 
M \to Gc$. Then $\tau'$ is pluriharmonic if and only if $f$ is a 
half isotropic ppmc immersion.
\endtheorem 

\proof
Recall from (6) that $d\tau' = \alpha|_{T^c \ox T'} \in \Hom(T^c \ox 
T', T'' + N)$. We compute $(Dd\tau')^{(1,1)}$. Let $X,Z$ be 
(1,0)-vector fields and $\bar Y$ a (0,1)-vector field. Then
$$
      (D_Xd\tau')(\bar Y).Z =  \pi''\d_X(\alpha(\bar Y,Z))  
                    + (D^N_X\alpha)(\bar Y,Z)
\leqno\ (7)
$$
where $\pi''$ is the projection onto $\tau'' \subset \C^n$ (which 
we identify with $T''$) and $D^N_X\alpha$ denotes the normal 
derivative of $\alpha$. Hence $\tau'$ is pluriharmonic if and only if both 
terms at right hand side vanish. The first term is zero if and only if $0 = 
\<\d_X(\alpha(\bar Y,Z)),W\> = -\<\alpha(\bar Y,Z),\alpha(X,W)\>$ 
for all $W \in T'$ which means that $\alpha(T',T') \in N^1 = 
(N^o)^\perp$. The vanishing of the second term is precisely the 
ppmc condition.
\endproof

\medskip\ni
{\bf Remark 1.}
It might seem more natural to use the embedding $j : Z_1 \subset 
Gc$ in order to prove the above theorem; clearly, $\tau' = j \o 
\tau_1 : M \to Gc$ is pluriharmonic if and only if $\tau_1$ is pluriharmonic and 
$(\tau_1^*\beta)^{(1,1)} = 0$ where $\beta$ denotes the second 
fundamental form of $Z_1 \subset Gc$. In fact, $(\tau_1^*\beta)(X,\bar Y)$ 
is given by the first summand at the right hand side of (7).  
Proving this involves computing the normal space and the second 
fundamental form of the submanifold $Z_1 \subset Gc$.

\smallskip
\ni
{\bf Remark 2.}
Half isotropic ppmc immersions are studied in [FT]. Such an 
immersion is always minimal in a sphere $S_r^{n-1}$ if it is 
substantial and indecomposable as a submanifold. In fact, the 
mean curvature vector $\eta = {1\over2m}\trace\alpha = {1\over2m}
\trace\alpha^{(1,1)} \in N^o$ is umbilic which can be seen as 
follows. First of all, $\eta$ is a parallel normal vector field 
since $\alpha^{(1,1)}$ is parallel. Further, the symmetric 
bilinear form $\alpha_\eta(x,y) = \<\alpha(x,y),\eta\>$ is 
parallel on $T' \ox T''$ and vanishes on $T' \ox T'$ and on $T'' \ox 
T''$ since $\alpha$ maps these bundles into $N^1$ which is 
perpendicular to $\eta$. Thus the corresponding Weingarten map 
$A_\eta$ is parallel. If $A_\eta$ 
had two different eigenvalues, the corresponding 
eigenspace distributions would give an extrinsic splitting of 
the immersion. Hence $A_\eta = \kappa\.I$ for some constant 
$\kappa > 0$. Therefore $m = f + {1\over\kappa}\,\eta$ is a 
constant point in $\R^n$, and $f(M)$ is contained in the sphere 
of radius ${1\over\kappa}$ centered at $m$. Since the mean 
curvature vector $\eta$ is normal to this sphere, the immersion
is minimal.

\bigskip
\bigskip
\leftline {\mittel 5. Isotropy}

\medskip
\noindent
We have seen that a ppmc immersion $f : M \to \R^n$ has an 
associated family of isometric immersions $f_\t$ with rotated 
second fundamental forms (cf. Equation (1)). It may happen that 
this family is trivial, i.e. $f_\t = f$ for all $\t$ (up to 
Euclidean motions) which implies some symmetry for the second 
fundamental form $\alpha$. In fact we see from (1) that $f_\t = 
f$ for all $\t$ if and only if there is a family of parallel vector bundle 
automorphisms $\psi_\t : N \to N$ with
$$
                  \psi_\t \o \alpha = \alpha_\t
\leqno\ (8)
$$
where $\alpha_\t(x,y) = \alpha(\r_\t x,\r_\t y)$ as before. We 
will call such an immersion {\it isotropic}. By the following 
theorem (cf. [ET2]), this property can be be read off from the 
{\it components of} $\alpha$:
$$\eqalign{
           \alpha^{(2,0)}(x,y) &= \alpha(\pi'x,\pi'y),      \cr
           \alpha^{(1,1)}(x,y) &= \alpha(\pi'x,\pi''y)
                                 + \alpha(\pi''x,\pi'y),    \cr
           \alpha^{(0,2)}(x,y) &= \alpha(\pi''x,\pi''y).    \cr
}$$

\theorem{Theorem 8.} An isometric K\"ahler immersion $f : M \to 
\R^n$ is isotropic ppmc if and only if there is a parallel 
orthogonal decomposition of the complexified normal bundle $N^c = 
N' \oplus N^o \oplus N''$ such that the parallel subbundles $N'$, 
$N^o$ and $N''$ contain the values of $\alpha^{(2,0)}$, 
$\alpha^{(1,1)}$ and $\alpha^{(0,2)}$, respectively. 
\endtheorem

\proof
If $f$ is isotropic ppmc, then the components of $\alpha$
take values in the eigenbundles of $\psi_\t$ corresponding to the 
eigenvalues $e^{2i\t}$, $1$ and $e^{-2i\t}$. They will be called 
$N'$, $N^o$ and $N''$. Since $\psi_\t$ is parallel, 
they form a parallel orthogonal decomposition of $N^c$. 
Vice versa, if such a decomposition of $N^c$ is given, we can 
define a parallel bundle automorphism $\psi_\t : N \to N$ by 
putting $\psi_\t = I$ on $N^o$ and $\psi_\t = e^{\pm 2i\t}I$ on 
$N'$ and $N''$, and we obtain Equation (8) which is equivalent to 
$f$ being isotropic ppmc.
\endproof

\remark
Theorem 8 implies in particular that isotropic ppmc immersions 
are half isotropic (cf. Ch. 4) since $\alpha^{(2,0)}$ takes 
values in in $N'$ which is perpendicular to $N^o$. Hence, by Remark 2
in Ch. 4 we may assume that an isotropic 
ppmc immersion takes values in a sphere $S^{n-1} \subset \R^n$. 
Thus Theorem 6 applies and in principle, we can obtain these 
immersions from their Gau\3 maps.

\medskip
By Theorem 3, isotropy of a ppmc immersion $f : M \to \R^n$ 
implies the isotropy of its Gau\3 map $\tau : M \to Gr$. The 
converse statement however cannot be true: If $f : M \to \R^n$ is 
{\it pluriminimal}, i.e. a pluriharmonic isometric immersion, its 
associated family $f_\t$ satisfies
$$
                       df_\t = df \o \r_\t
$$
up to a rigid motion of $\R^n$ (cf. [ET2]), hence we also 
conclude $\tau_\t = \tau$ (another argument for the isotropy of 
$\tau$ will be given below). But we will see in the next theorem 
that these are essentially the only two cases where the Gau\3 map 
is isotropic.

We need some preparations. For any complex vector bundle $E 
\subset M \x \C^n$, let us define a linear map $d : T^c \to 
\Hom(E,E^\perp)$ (the {\it differential} or {\it shape operator} of 
$E$) by assigning to each vector $v \in T^c$ and any section $s$ 
of $E$ the $E^\perp$-component of $\d_vs$. According to the 
splitting $T^c = T' + T''$, the differential splits as $d = d' + 
d''$.

\theorem{Lemma 4.} For any isotropic ppmc immersion $f : M \to 
\R^n$ we get the following chain of differentials:
$$\eqalign{
  d' &: N'' \to \tau'' \to N^o \to \tau' \to N' \to 0     \cr
 d'' &: N' \to \tau' \to N^o \to \tau'' \to N'' \to 0     \cr
}$$
\endtheorem

\proof
Since $N', N^o, N''$ are parallel subbundles of $N^c$, being 
eigenbundles of the parallel bundle automorphism $\psi_\t : N \to 
N$, the differential of any of them takes values in $\tau^c$. 
Similarly, $\tau'$ and $\tau''$ are mapped into $N^c$, being 
parallel subbundles of $\tau^c$. Hence $d'\tau'' = \alpha(T',T'') 
= N^o$. Further, $\<d'N'',\tau''\> = \<N'',d'\tau''\> = 
\<N'',N^o\> = 0$ and consequently $d'N'' \subset \tau''$ since 
$\tau'' \subset \tau^c$ is maximal isotropic. Next, 
$\<d'N^o,\tau'\> = \<N^o,d'\tau'\> = \<N^o,N'\> = 0$, thus $d'N^o 
\subset \tau'$. Further, $d'\tau' = N'$. Finally, $\<d'N',\tau'\> 
= \<N',N'\> = 0$ since $N'$ is isotropic (being perpendicular to 
$N'' = \overline{N'}$), and $\<d'N',\tau''\> = \<N',N^o\> = 0$, 
thus we get $d'N' = 0$. This proves the first chain of 
differentials. The second one follows by complex conjugation.
\endproof

\theorem{Lemma 5.} Let $M = M_1 \x M_2$ be a Riemannian product 
of K\"ahler manifolds and $f : M \to \R^n$ an isometric 
immersion. Let $x_1 \in TM_1$ and $x_2 \in TM_2$. Then 
$|\alpha^{(1,1)}(x_1,x_2)| = |\alpha^{(2,0)}(x_1,x_2)|$. In 
particular $\alpha(x_1,x_2) = 0$ if and only if $\alpha^{(2,0)}(x_1,x_2) = 
0$. If this holds for all such $x_1,x_2$, the splitting is extrinsic.
\endtheorem

\proof
Since all mixed curvature tensor components of the Riemannian product
$M$ are zero, we obtain from the Gau\3 equation for any 
$y_1 \in T^cM_1$ and $y_2 \in T^cM_2$
$$\eqalign{ 
               0 = \<R(y_1,\bar y_1)y_2,\bar y_2\> 
         &= \<\alpha(y_1,\bar y_2),\alpha(\bar y_1,y_2)\>
         - \<\alpha(y_1,y_2),\alpha(\bar y_1,\bar y_2)\>  \cr
         &= |\alpha(y_1,\bar y_2)|^2 - |\alpha(y_1,y_2)|^2.\cr
}$$
Thus $|\alpha(y_1,y_2)| = |\alpha(y_1,\bar y_2)|$ and in 
particular, putting $y_1 = \pi'x_1$ and $y_2 = \pi'x_2$, we get 
$$
     |\alpha(\pi'x_1,\pi'x_2)| = |\alpha(\pi'x_1,\pi''x_2)|.
$$
The extrinsic splitting is obvious if $\alpha(TM_1,TM_2) = 0$.
\endproof

\theorem{Lemma 6.} Let $H \subset O(2m)$ be a group acting on $V 
= \R^{2m}$ and let $J,\tilde J \in O(2m)$ be two $H$-invariant 
complex structures on $V$. Then there is an $H$-invariant 
decomposition $V = \sum_j V_j$ such that on each $V_j$ we have 
either $\tilde J = \pm J$ or there is an $H$-invariant 
quaternionic structure on $V_j$.
\endtheorem

\proof
Using the complex structure $J$, we consider $\R^{2m}$ as a 
complex vector space, and we decompose $\tilde J$ into its 
complex linear and antilinear components (called $L$ and $A$). 
Hence $\tilde J = L + A$ with $L = {1\over2}(\tilde J - J\tilde J 
J)$ and $A = {1\over2}(\tilde J + J\tilde J J)$. From $\tilde J^2 
= -I$ we get
$$
                    -I = L^2 + A^2 + LA + AL,
$$
and since $L^2 + A^2$ is linear while $LA + AL$ is antilinear, 
this implies $L^2 + A^2 = -I$ and $LA + AL = 0$. Since both $L$ 
and $A$ are antisymmetric, $L^2$ and $A^2 = -L^2 - I$ are 
symmetric and decompose $V = \R^{2m}$ into common eigenspaces 
$W_1,...,W_r$ with non-positive real eigenvalues. Let $W = W_j$ 
be any of these eigenspaces and $-c^2,-s^2$ with $c^2 + s^2 = 1$ 
the corresponding eigenvalues of $L^2$ and $A^2$. If $s = 0$, we 
have $A = 0$ and $\tilde JJ = J\tilde J$ on $W$. Thus there is an 
$H$-invariant splitting $W = W_+ + W_-$ with $\tilde J = J$ on 
$W_+$ and $\tilde J = -J$ on $W_-$. If $s \neq 0$, we may put 
$J_2 = {1\over s}A$ and obtain $(J_2)^2 = {1\over s^2}A^2 = -I$. 
This is an antisymmetric complex structure, hence orthogonal 
(since $(J_2)^T = -J_2$ and $(J_2)^2 = -I$ imply $(J_2)^TJ_2 = 
I$), and $J_2$ anti-commutes with $J$. Thus $J_1 := J$ together 
with $J_2$ and $J_3 := J_1J_2$ form an $H$-invariant quaternionic 
structure on $W$.
\endproof

\theorem{Corollary.} Let $M$ be a locally irreducible Riemannian
manifold with two linear independent parallel almost complex
structures. Then $M$ is locally {\it hyper-K\"ahler}, i.e. 
locally there exist three anti-commuting parallel almost complex 
structures on $M$.
\endtheorem

\proof
We apply Lemma 6 for $V = T_pM$ where $H$ is the local holonomy 
group of $M$ at the point $p$. By assumption this acts 
irreducibly, so the $H$-invariant decomposition $V = \sum V_j$ 
must be trivial. Since the two almost complex structures are 
linearly independent, we get a quaternionic structure 
$(J_1,J_2,J_3)$ on $T_pM$ which is invariant under the local 
holonomy group and thus allows a parallel extension on a 
neighborhood of $p$.
\endproof

\theorem{Theorem 9.} Let $M$ be a K\"ahler manifold such that no 
local factor of $M$ is  hyper-K\"ahler, and let $f : M \to 
\R^n$ be an isometric immersion with Gau\3 map $\tau : M \to Gr$. 
Then $\tau$ is isotropic pluriharmonic if and only if $f$ is 
either pluriminimal or isotropic ppmc.
\endtheorem

\proof
The map $\tau : M \to Gr$ is isotropic pluriharmonic if and only if there is 
a holomorphic superhorizontal lift $\hat\tau : M \to Z$ into 
some flag manifold $Z$ fibering over $Gr$ (cf. [ET2]). We 
classify these flag manifolds in the appendix and obtain $Z = 
Z_r$ for some $r \in \N$, where $Z_r$ is the set of all 
$(2r+1)$-tuples of complex subspaces $E_{-r},...,E_r$ with
given dimensions forming 
an orthogonal decomposition $\C^n = \sum\limits_{j=-r}^r E_j$ 
such that $E_{-j} = \overline{E_j}$ for all $j$. Thus the lift 
$\hat\tau$ is a ``moving" orthogonal decomposition $(E_{-
r},...,E_r)$ of subbundles $E_j \subset M \x \C^n$ with $E_{-j} = 
\overline{E_j}$, and the fact that $\hat\tau$ is holomorphic 
superhorizontal means that $d'E_j = E_{j+1}$. Since $\hat\tau$ is 
a lift of $\tau$, we have either $\tau^c = E_{{\rm even}}$ or 
$\tau^c = E_{{\rm odd}}$ where $E_{{\rm even}} = \sum\limits_{j+r\ 
{\rm even}} E_j$ and $E_{{\rm odd}} = \sum\limits_{j+r\ {\rm odd}} 
E_j$.

Now $f : M \to \R^n$ is pluriminimal if and only if $\tau'$ is holomorphic 
which means $d''\tau' = 0$. Consequently $d''N \subset \tau'$ 
(since $\<d''N,\tau'\> = \<N,d''\tau'\> = 0$) and $d''\tau'' 
\subset N$ (since $\tau'' \subset \tau^c$ is parallel), hence
$$
            d'' : \tau'' \to N \to \tau' \to 0, \ \ \ 
               d' : \tau' \to N \to \tau'' \to 0.
$$
Thus $\hat\tau = (\tau',N,\tau'')$ is a (super-)horizontal 
holomorphic lift into the corresponding flag manifold $Z_1$ (and 
in particular, $\tau$ is isotropic pluriharmonic).

If $f : M \to \R^n$ is isotropic ppmc, then 
$\hat\tau = (N'',\tau'',N^o,\tau',N')$ is a superhorizontal holomorphic lift 
into the corresponding flag manifold $Z_2$ (cf. Lemma 4).

\medskip
Conversely, let $f : M \to \R^n$ be any K\"ahler immersion such that the 
Gau\3 map $\tau : M \to Gr$ is isotropic pluriharmonic and let 
$\hat\tau = (E_{-r},...,E_r)$ be the holomorphic superhorizontal 
lift of $\tau$. Then $\tau^c = E_{-r'} + E_{-r'+2} + ... + 
E_{r'}$ where $r' \in \{r-1,r\}$, and since $d'E_j \subset 
E_{j+1}$ and $d''E_j \subset E_{j-1}$, the subbundles $E_j$ of 
$\tau^c$ are parallel. Let  $\tau_j^c = E_j+E_{-j}$. Then 
$\tau_j^c = \tau_j \ox \C$ for some parallel real subbundle 
$\tau_j \subset \tau$, and $E_{\pm j} = (I \mp iJ_j)\tau_j$ for a 
parallel complex structure $J_j$ on $\tau_j$, if $j \neq 0$. By 
the Corollary of Lemma 6 and the present assumption we may assume 
$J_j = \pm J$ where $J$ is the complex structure of $TM$, 
transplanted by $df$ onto $\tau$. (Maybe we yet have to split 
$\tau_j$ into holonomy irreducible subbundles.) Thus $E_j = 
\pi'(\tau_j)$ or $E_j = \pi''(\tau_j)$. If $E_i = \tau_i'$ and 
$E_j = \tau_j'$ for some $i \neq \pm j$, using the symmetry of 
$\alpha$ we have
$$
\alpha(E_i,E_j) = d(E_i).E_j \subset E_{j-1} \cap E_{i-1} = 0
$$
and likewise, if $E_i = \tau_i'$ and $E_j = \tau_j''$, we have
$$
\alpha(E_i,E_{-j}) = d(E_i).E_{-j} \subset E_{-j-1} \cap 
                          E_{i-1} = 0.
$$
In both cases we get $\alpha^{(2,0)}(\tau_i,\tau_j) = 0$ which by 
Lemma 5 is equivalent to $\alpha(\tau_i,\tau_j) = 0$ (recall that 
the parallel subbundles $\tau_j$ define a local Riemannian 
product structure on $M$). So we see that the splitting is also 
extrinsic and we may assume $r' = 1$.

\smallskip\noindent
The remaining possibilities for our 
moving flag are the following four cases: $(\tau',N,\tau'')$, 
$(\tau'',N,\tau')$, $(N'',\tau',N^o,\tau'',N')$, and 
$(N'',\tau'',N^o,\tau',N')$ (the bundles $N'$ and $N''$ are 
interchangeable). In the first case we have $d''\tau' = 0$, so 
$\tau'$ is holomorphic and hence $f$ is pluriminimal by Theorem 
4. The second case is equivalent to $d'\tau' = 0$ which means 
$\alpha^{(2,0)} = 0$. This implies $(D\alpha)^{(2,1)} = 0$ and 
hence $D\alpha = 0$ by the Codazzi equation, and in particular 
$f$ is a ppmc immersion. In fact these are the standard 
embeddings of compact Hermitian symmetric spaces (cf. Ch. 7). In 
the third case we get $d'\tau'' \subset N'$ and $d''\tau' \subset 
N''$. Hence $\alpha(T',T'') \in N' \cap N'' = 0$ and thus 
$\alpha^{(1,1)} = 0$. So we are back to the first case. Finally 
in the last case, $\alpha^{(2,0)}$, $\alpha^{(1,1)}$ and 
$\alpha^{(0,2)}$ take values in the parallel subbundles $N'$, 
$N^o$ and $N''$ which shows isotropy by Theorem 8.
\endproof

\bigskip
\bigskip
\leftline {\mittel 6. Isotropy and complex Gau\3 map}

\medskip
\noindent
Using the complex Gau\3 map with values in the complex 
Grassmannian, we can characterize isotropy avoiding the 
unpleasant extra condition of Theorem 9:

\theorem{Theorem 10.} A K\"ahler immersion $f : M \to \R^n$ is 
isotropic ppmc if and only if its complex Gau\3 map $\tau' : M 
\to Gc$ is isotropic pluriharmonic, but not holomorphic.
\endtheorem

\proof
Assume first that $f : M \to \R^n$ is isotropic ppmc. Then we 
have an orthogonal decomposition (``moving flag") $\C^n = N' 
\oplus \tau' \oplus Q$ with $Q := N^o + \tau'' + N''$ (where 
$\C^n$ denotes the trivial vector bundle $M \x \C^n$), and by 
Lemma 4 we have the differentials $d' : Q \to \tau' \to N' \to 0$ 
and $d'' : N' \to \tau' \to Q \to 0$. Thus the map $(Q,\tau',N')$ 
into the corresponding flag manifold over $Gc$ with the 
projection $(Q,\tau',N') \mapsto \tau'$ is horizontal and 
harmonic and thus $\tau'$ is isotropic pluriharmonic.

Conversely, let us assume that $\tau'$ is isotropic 
pluriharmonic, i.e. there is a one parameter group $\phi_\t \in 
\Aut(\tau'^*(TGc))$ with $\phi_\t \o d\tau' = d\tau' \o {\cal 
R}_\t$. By [ET2] we have a horizontal holomorphic lift 
$\hat\tau'$ of $\tau'$ into some flag manifold $Z$ over $Gc$, 
i.e. (cf. Appendix) there are decompositions $\tau' = \tau_1' 
\oplus ... \oplus \tau_r'$ and $(\tau')^\perp = P_1 \oplus ... 
\oplus P_{r+1}$ (where $P_1$ and $P_{r+1}$ might be zero) such 
that $d' : P_i \to \tau_i \to P_{i+1}$ and $d'' : P_{i+1} \to 
\tau_i \to P_i$ for $i = 1,...,r$. By the following argument we 
may assume $r = 1$ and thus $\hat\tau'$ is a ``moving 
decomposition" of the type
$$
               \C^n = P_1 \oplus \tau' \oplus P_2.
$$
In fact, the parallel decomposition $\tau' = \tau_1' \oplus ... 
\oplus \tau_r'$ induces a corresponding real parallel decomposition 
$TM = T_1 \oplus ... \oplus T_r$ and hence the manifold $M$ 
can be (locally) decomposed as a Riemannian product of 
K\"ahler manifolds. This splitting is even extrinsic: 
For any $x_i' \in T_i'$ and $x_j' \in T_j'$ we have (using $df$ 
to identify $T'$ and $\tau'$)
$$\eqalign{
\alpha(x_i',x_j') 
     &= d\tau'(x_i').x_j' \subset d'\tau_j \subset P_{j+1},    \cr
\alpha(x_j',x_i')
     &= d\tau'(x_j').x_i' \subset d'\tau_i \subset P_{i+1},    \cr    
}$$
hence $\alpha(x_i',x_j') = 0$ and thus $\alpha^{(2,0)}(T_i,T_j) = 
0$. But by Lemma 5 this implies $\alpha(T_i,T_j) = 0$. Hence we 
may assume $r = 1$. 

\medskip
Now we claim for any (1,0) vector fields $X,Y,Z$ (while still 
identifying $TM$ with $\tau$)
$$
     (D_{\Z}d\tau').\X.Y 
     = (D^N_{\Z}\alpha)(\X,Y) + (\d_{\Z}(\alpha(\X,Y))^{T''}
\leqno\ (8)
$$
In fact, recall from (6) (Lemma 2) that
$$
      d\tau' : TM \to \tau'^*(TGc) = \Hom(\tau',\tau'^\perp), \ \
            d\tau'.V.Y = (\d_{V}Y)^{T''+N} = \alpha(V,Y)
$$
for any $V \in T^c$ and $Y \in T'$. Then 
$$
(D_{\Z}d\tau').\X.Y = 
(\d_{\Z}(d\tau'.\X.Y))^{T''+N} - d\tau'.D_{\Z}X.Y - d\tau'.X.D_{\Z}Y
\leqno\ (9)
$$
Now we may replace $d\tau'$ by $\alpha$. Consider the right hand 
side of (9) (``rhs(9)"). The first term splits into its 
components with respect to $T''$ and $N$. Its $N$-component 
together with the 2nd and 3rd terms gives 
$(D^N_{\Z}\alpha)(\X,Y)$ (which is the first term of rhs(8)) 
while the remaining term $\d_{\Z}(\alpha(\X,Y)^{T''}$ is the 
second summand of rhs(8). Thus Equation (8) is proved.

On the other hand we have
$$
      (D_{\Z}d\tau').\X = D_{\Z}(d\tau'.\X) - d\tau'.D_{\Z}\X
\leqno\ (10)
$$
If $\tau'$ is isotropic pluriharmonic, then both terms at rhs(10) 
are eigenvectors of $\phi_\theta$ with respect to the eigenvalue 
$e^{-i\theta}$: the second one because $D_{\Z}\X \in T''$ and 
$\phi_\theta \o d\tau' = d\tau' \o \r_\theta$, and the first one 
because the eigenbundle of $\phi_\theta$ is parallel. Thus these 
vectors lift to (0,1) superhorizontal tangent vectors of $Z$ (cf. 
[ET2]) which map $P_2 \to T' \to P_1$. 

It follows that $(D_{\Z}d\tau').\X$ maps $T'$ into $P_1$, and 
since the first term of rhs(8) vanishes by the ppmc property, we 
conclude from (8) that 
$$
               (\d_{\Z}\alpha(\X,Y))^{T''} \in P_1
$$ 
Thus putting $T''_0 = T'' \cap P_1$ and letting $T''_1$ be the 
orthogonal complement of $T''_0$ in $T''$, we have 
$(\d_{\Z}\alpha(\X,Y))_{T''_1} = 0$, and therefore we obtain for 
all $W \in T'$ with $\bar W \in T''_1$:
$$
  \<\alpha(\X,Y),\alpha(\Z,W)\> = \<\d_{\Z}\alpha(\X,Y),W\> = 0.
$$
In other words, $\alpha(\Z,W) = 0$ for all $Z \in T'$ which says 
that $W$ and hence all of $T''_1$ lies in the subbundle
$$
\ker \alpha^{(1,1)} := 
     \{W \in T';\ \alpha(\Z,W) = 0\ \ \forall Z \in T'\}.
$$
By parallelity of $\alpha^{(1,1)}$, this is a parallel subbundle 
of $T'$ which can be split off, using Lemma 5 (yielding a 
pluriminimal factor). Thus we may assume that $\ker\alpha^{(1,1)} 
= 0$ and hence $T''_1 = 0$, i.e. $T'' \subset P_1$. 

Similar as in (10) we have that $(D_{Z}d\tau').X$ is in the 
$e^{i\theta}$-eigenspace of $\phi_\t$ whose elements map $T'$ 
into $P_2$, and as in (8) we have
$$
     (D_{Z}d\tau').X.Y 
          = (D^N_{Z}\alpha)(X,Y) 
            + (\d_{Z}(\alpha(X,Y))^{T''} \in P_2.
\leqno\ (11)
$$
But the second term of rhs(11) is in $T'' \subset P_1$ while the 
first one is in $N \perp T''$. Hence the sum can be perpendicular 
to $T''$ (recall that $P_2 \perp P_1 \supset T''$) only if its 
$T''$-component (the second term of rhs(11)) vanishes. Taking the 
inner product of this term with any $W \in T'$ we obtain
$$
                \<\alpha(X,Y),\alpha(Z,W)\> = 0
$$
for arbitrary $X,Y,Z,W \in T'$. Thus $\<N',N'\> = 0$ or in other 
words $N' \perp \overline{N'} =N''$. Since $f$ is already half 
isotropic (cf. Theorem 7), we also have $N' \perp N^o$. Now the 
proof is finished by Theorem 8.
\endproof

\bigskip
\bigskip
\goodbreak
\leftline{\mittel 7. Examples.}

\medskip
\noindent
Clearly, if $f_i : M_i \to \R^{n_i}$ are any two ppmc K\"ahler 
immersions ($i = 1,2$), then so is $f = f_1\x f_2 : M_1 \x M_2 
\to \R^{n_1+n_2}$. Therefore it is enough to study ppmc 
immersions $f : M \to \R^n$ which are {\it irreducible}, i.e. 
they do not split as above, and {\it substantial}, i.e. their 
image is not contained in any proper affine subspace of $\R^n$. 
Three classes of such immersions are known:

\smallskip
\item{(1)}  surfaces with nonzero parallel mean curvature vector,
\item{(2)}  pluriminimal submanifolds,
\item{(3)}  extrinsic symmetric K\"ahler immersions.

\smallskip
\noindent
Class (1) has been investigated by S.T.Yau [Y]; these examples 
occur only in $\R^3$ or $S^3$ unless they are minimal surfaces in 
$S^{n-1}$. Class (2) contains many examples in all dimensions, 
cf. [DG] and the recent paper [APS]. We will now briefly describe 
Class (3).      

\medskip
Recall that an isometric (irreducible, substantial) immersion $f 
: M \to \R^n$ is called {\it extrinsic symmetric} if the full 
second fundamental form $\alpha \in \Hom(TM \otimes TM, N)$ is 
parallel. These immersions have been classified by D.Ferus ([F], 
also cf. [EH]). It is not difficult to see that $\alpha$ is 
parallel if and only if $f$ is invariant under reflection at each of its 
normal spaces. In particular all point reflections or geodesic 
symmetries on $M$ extend to (extrinsic) isometries, hence $M$ is 
globally symmetric. Moreover, $M$ is isotropy irreducible, i.e. 
the full extrinsic isotropy group of $M$ acts irreducibly on the 
tangent space (cf. [EH]). The corresponding Gau\3 map $\tau : M 
\to Gr$ is a totally geodesic isometric immersion of the 
symmetric space $M$ into the real Grassmannian $Gr$. In fact, 
since $\tau$ is equivariant and $M$ is isotropy irreducible, it 
is isometric (up to a scaling factor). Moreover, the image of 
$\tau$ is invariant under the corresponding point reflections 
of $Gr$ and thus totally geodesic; note that the point reflection of 
the Grassmannian at some $\tau(p) \in Gr$ is just the reflection at 
the normal space $\tau(p)^\perp = N_p$.

Hence, if $f : M \to \R^n$ is an extrinsic symmetric immersion 
which is also K\"ahler (with almost complex structure $J$), then 
$f$ is clearly ppmc since the parallelity of $\alpha^{(1,1)}$ is 
a weaker condition. Moreover, if $f$ is also substantial and 
irreducible, it is isotropic. To see this recall that a symmetric 
space $M$ with a K\"ahler metric is in fact {\it Hermitian 
symmetric}, i.e. the rotations $\r_\t(p) = \cos(\t)I + \sin(\t)J$ 
on $T_pM$ for any $p \in M$ extend to isometries $\rho_\t$ on $M$ 
fixing $p$. But these isometries are generated by point 
reflections which extend to orthogonal linear maps on $\R^n$, 
hence $\rho_\t$ also extends to some $A_\t \in O(n)$ with $f \o 
\rho_\t = A_\t \o f$. We put $\psi_\t(p) = A_{2\t}|_{N_p}$. Since 
$A_\t$ (being an extrinsic isometry) commutes with $\alpha$, we 
obtain
$$
         \psi_\t(\alpha(v,w)) = \alpha(\r_\t v,\r_\t w) 
\leqno\ (8)
$$
for all $v,w \in T_pM$. In particular this equation implies that 
$p \mapsto \psi_\t(p)$ is parallel (as an endomorphism of the 
normal bundle $N$), since so are $\r_\t$ and $\alpha$ and since
$N = \alpha(TM \otimes TM)$. Thus $f$ is isotropic ppmc.

Since $\psi_{\pi} = I$ by (8), the eigenvalues of $\psi_{\pi/2}$ 
can only be $\pm 1$. Accordingly, class (3) has two subclasses: 
If $1$ is the only eigenvalue, i.e. $\psi_{\pi/2} = I$, then we 
get from (8)
$$
                   \alpha(Jv,Jw) = \alpha(v,w)
$$
for all $v,w$, hence $\alpha^{(2,0)} = 0$. These immersions have 
been characterized already by Ferus [F]: They are the so called 
{\it standard embeddings} of an Hermitian symmetric space $M = 
G/K$ into the Lie algebra $\g$ of $G$ via the map $p \mapsto J_p$
(recall that the complex structure $J_p$ on $T_pM$ is a skew-symmetric 
derivation of the curvature tensor of $M$ at $p$, hence it 
extends to an infinitesimal isometry, i.e. to an element of 
$\g$).

In the remaining examples, the eigenvalue $-1$ occurs for $\psi_{\pi/2}$. 
Inspection shows that these are precisely the extrinsic 
symmetric $2:1$ immersions of $Gr_2^+ = G_2^+(\R^N)$, the Grassmannian 
of {\it oriented} 2-planes in $\R^N$, factorizing over the 
ordinary real Grassmannian $Gr_2$. In fact, $Gr_2^+$ is an Hermitian 
symmetric space (which can be identified with the complex quadric 
$\{[z] \in \C P^{N-1};\ \<z,z\> = 0\}$ via the map $E = 
\Span\{x,y\} \mapsto [x+iy]$, where $(x,y)$ is any oriented 
orthonormal basis of the oriented plane $E \subset \R^N$). We put 
$f = \tilde f \o \pi$ where $\pi : Gr_2^+ \to Gr_2$ is the 
canonical projection and $\tilde f : Gr_2 \to S(N)$ the usual 
(extrinsic symmetric) embedding of the Grassmannian into the 
space of symmetric real $N \x N$-matrices by assigning to each 
plane $E \in Gr_2$ the orthogonal projection of $\R^N$ onto $E$. 
In this case, the $(-1)$-eigenspace is $2$-dimensional. The 
easiest example is the Veronese immersion 
$$
                 S^2 \to \R P^2 \to S^4 \subset 
          \R^5 \cong \{X \in S(3);\ \trace X = 1\}.
$$

It is an open problem how to construct further classes of examples. 
Using our Theorem 6, we hope that a better understanding of 
horizontal pluriharmonic maps into $Z_1$ will lead to new ppmc 
immersions.

\bigskip
\bigskip
\leftline {\mittel Appendix: Canonical embeddings of flag 
manifolds}

\medskip
\noindent
Let $G$ be a compact Lie group with Lie algebra $\g$, and let 
$\g^c = \g \ox \C$ be the complexification of $\g$. We consider adjoint 
orbits (``flag manifolds'') $Z = Ad(G)\xi$ for $\xi \in \g$. 
An orbit can always be represented as a coset space $G/H$ 
where $H$ is the stabilizer subgroup; in the present case 
$H = C(\xi) =\{g \in G;\ Ad(g)\xi = \xi\}$ is 
the centralizer of $\xi$. More precisely, $Z$ is the image of the
equivariant embedding $j_\xi : G/H \to \g$, $j_\xi(gH) = Ad(g)\xi$.
Of course, if we fix $H$, many $\xi \in \g$ may have $H$ as centralizer and 
give different embeddings $j_\xi$ of the same coset space $G/H$, 
but there are distinguished such $\xi$: We call $\xi \in \g$ a 
{\it canonical element} and $j_\xi$ a {\it canonical embedding} 
of $G/H$ for $H = C(\xi)$ if

\smallskip
\item\item{{\bf C1}} \ \ The eigenvalues of ${1\over i}ad(\xi)$ are 
            integers (where $i = \sqrt{-1}$),
\item\item{{\bf C2}} \ \ $\g_1 + \g_{-1}$ generates $\g^c$, where $\g_k 
            \subset \g^c$ denotes the $k$-eigenspace of ${1\over 
            i}ad(\xi)$. \footnote{$^*)$}
     {A canonical element $\xi$ is not uniquely determined by $H$. 
     But there is only one such $\xi$ (up to adding an element in 
     the center of $\g$) in any Weyl chamber $C$ of $\g$ which is
     adjacent to the subtorus $T'$ centralized by $H$ 
     (where ``adjacent'' means that $\bar C \cap \T'$ contains an
     open subset of $\T' = L(T')$). In fact   
     $\xi = \sum_{j \in J} \alpha_j^*$, where $\{\alpha_1,..., 
     \alpha_l\}$ are the simple roots of $\g$ corresponding to
     $C$ and $\alpha_1^*,..., \alpha_l^*$ the
     dual root vectors (i.e. 
     $\alpha_j(\alpha_k^*) = \delta_{jk}$) and where $J = \{j \in 
     \{1,...,l\};\ \g_{\alpha_j} \cap \h = 0\}$ (cf. [BR], p.42).
     Using this extra structure we can represent $G/H$ as the complex
     coset space $G^c/P$ for the parabolic subgroup $P = \{g \in G^c;\
     Ad(g)\xi \in \xi + \sum_{k > 0} \g_k\}$, and our definition of
     ``canonical element'' agrees with that of [BR], p.41.}
\smallskip
\ni
The Jacobi identity implies $[\g_j,\g_k] \subset \g_{j+k}$. Since 
$\g_1+\g_{-1}$ is a generating subspace and $\g_{-j} = 
\overline{\g_j}$, the eigenvalues of ${1\over i}ad(\xi)$ form a set 
$\{-r,...,r\}$ for some positive integer $r$ (called the {\it height} 
of the flag manifold) where $\g_0 = \h^c$ is the complexified Lie 
algebra of $H$, and we have a direct decomposition
$
                   \g^c = \sum\limits_{j=-r}^r \g_j .
$

The flag manifold $Z = G/H$ fibres over a symmetric space $S = 
G/K$ defined by the corresponding (complexified) Cartan 
decomposition as follows:
$$
              \k^c = \sum_{j\ {\rm even}} \g_j,\ \ 
                \p^c = \sum_{j\ {\rm odd}} \g_j .
\leqno\ (A1)
$$
In fact, the Cartan relations $[\k,\k] \subset \k$, $[\k,\p] 
\subset \p$, $[\p,\p] \subset \k$ are obvious from $[\g_j,\g_k] 
\subset \g_{j+k}$, and clearly $\h^c = \g_0 \subset \k^c$. Thus 
$Z$ defines a unique symmetric space $S$ which is {\it inner}, 
i.e. its symmetry is an inner automorphism (namely 
$Ad(e^{\pi\xi})$). But conversely there are {\it several}
flag manifolds which fibre 
over $S$ as described. As an example we shall determine all canonical 
elements and corresponding flag manifolds over complex and  
real Grassmannians, using only elementary linear algebra.

\medskip

First let $G = U_n$ the unitary group. Then $\g = \u_n$ is the 
space of skew-Hermitian matrices. Any $\xi \in \g$ determines an 
orthogonal eigenspace decomposition of $\C^n$, and the 
eigenvalues are imaginary. Thus there is an orthogonal 
decomposition $\C^n = \sum_{j=1}^m E_j$ such that $\xi = i\. 
\sum_{j=1}^m \lambda_j E_j$ for real numbers $\lambda_1 < 
\lambda_2 < ... < \lambda_m$, where for any subspace $E \subset 
\C^n$ we use the same symbol $E$ to denotes the orthogonal 
projection matrix onto $E$. If $E,F \subset \C^n$ are subspaces 
with $E \perp F$, we embed $\Hom(E,F)$ into $\End(\C^n) = \g^c$ 
by putting $L|_{E^\perp} = 0$ for any $L \in \Hom(E,F)$. Then we 
have for any $L_{EF} \in \Hom(E,F)$:
$$
         [E,L_{EF}] = -L_{EF},\ \ \ [F,L_{EF}] = L_{EF}.
\leqno\ (A2)
$$
Thus for all $L_{jk} \in H_{jk} := \Hom(E_j,E_k)$ we obtain
$$
         ad(\xi)L_{jk} = i\.(\lambda_k-\lambda_j)\.L_{jk}.
\leqno\ (A3)
$$
Hence, if $\xi$ is canonical, then $\lambda_k-\lambda_j$ are 
integers for all $j,k$, by Property C1. Next we claim 
$\lambda_{j+1} -\lambda_j = 1$ for all $j$. This is due to 
Property C2 saying that $\g_1 + \g_{-1}$ generates $\g^c$. In 
fact, if $\lambda_{k+1} - \lambda_k \geq 2$ for some $k$, we may 
decompose $\C^n = E \oplus F$ with $E = \sum_{j=1}^k E_j$ and $F 
= \sum_{l=k+1}^m E_l$. Then $\lambda_l-\lambda_j \geq 2$ for all 
$j \in \{1,...,k\}$ and $l \in \{k+1,...,m\}$, and hence 
$H_{jl} = \Hom(E_j,E_l)$ and $H_{lj} = \Hom(E_l,E_j)$ belong to 
some $\g_k$ with $|k| \geq 2$. In other words, $\g_1 + \g_{-1}$ 
is contained in $\Hom(E,E) \oplus \Hom(F,F)$ which is a proper Lie 
subalgebra of $\u_n^c$. This contradicts Property C2. Thus we 
have seen (the converse statement is obvious):

\theorem{Proposition A1.} An element $\xi \in \g = \u_n$ is 
canonical if and only if $\xi = i(\lambda_0\.I + \sum_{j=1}^m j\.E_j)$ for 
some orthogonal decomposition $\C^n = \sum_{j=1}^m E_j$ and any 
$\lambda_0 \in \R$. Then $\g_k = \sum_j H_{j,j+k}$.
\endtheorem

The corresponding flag manifold is a ``classical" flag manifold 
$Z$ consisting of all orthogonal decompositions of $\C^n$ with the 
same dimensions as $E_1,...,E_r$, and $Z$ is embedded as the adjoint 
orbit $Ad(U_n)\xi$. What is the corresponding symmetric space $S$ 
over which $Z$ fibres? Let us put $E_{odd} = \sum_{j\ {\bf odd}} 
E_j$ and $E_{ev} = \sum_{j\ {\bf even}} E_j$. Then we have
$$
        \k^c = \End(E_{ev}) \oplus \End(E_{odd}),\ \ \ \
    \p^c = \Hom(E_{ev},E_{odd}) \oplus \Hom(E_{odd},E_{ev}).
\leqno\ (A4)
$$
This is the complexified Cartan decomposition of a symmetric 
space, namely the Grassmannian of all subspaces in $\C^n$ with 
the same dimension as $E_{ev}$ (or as $E_{odd})$. 

\bigskip
Now let $G = SO_n$ be the orthogonal group which we consider as a 
subgroup of $U_n$. Let $\xi \in \so_n \subset \u_n$. As before, 
we have $\xi = i\.\sum_{j=1}^m \lambda_j E_j$ for some orthogonal 
decomposition $\C^n = \sum_j E_j$  where $\lambda_1 < ... < 
\lambda_m$ are real. But now $\xi$ is a real matrix, i.e. we also 
have $\xi = \bar\xi = -i\.\sum_j \lambda_j \overline{E_j}$. Since 
the projections $E_j$ are linearly independent and nonnegative, 
there is a permutation $\sigma$ of $\{1,...,m\}$ such 
that $\overline{E_j} = E_{\sigma j}$ and $\lambda_{\sigma j} = -
\lambda_j$. Thus 
$$
   \hat H_{jk} := \Hom(E_j,E_k) + \Hom(E_{\sigma k},E_{\sigma j})
$$ 
is the eigenspace of $ad(\xi)$ corresponding to 
the eigenvalue $\lambda_k-\lambda_j$, according to $(A3)$. Now 
$\so_n^c = \{A \in \C^{n\x n};\ A^T = -A\}$ is generated as a 
vector space by $M_{jk}:= L_{jk} - (L_{jk})^T$ for all $L_{jk} 
\in \Hom(E_j,E_k)$ and all $j,k \in \{1,...,m\}$. We claim that 
$M_{jk} \in \hat H_{jk}$.

In fact, it is sufficient to show that $(L_{jk})^T \in 
\Hom(E_{\sigma k},E_{\sigma j})$. Put $y = (L_{jk})^Tx$ for some $x 
\in \C^n$. Let us denote the symmetric inner product on $\C^n$ by 
$\<v,w\> = \sum v_jw_j$. Then for all $w \in \C^n$ we have 
$\<y,w\> = \<x,L_{jk}w\>$, and the latter is nonzero only if $w 
\in E_j$ and $x \in \overline{E_k}$. Moreover $\<y,w\> \neq 0$ 
implies $y \in \overline{E_j}$. Thus $(L_{jk})^T$ maps 
$\overline{E_k} = E_{\sigma k}$ into $\overline{E_j} = E_{\sigma j}$ 
and vanishes on the orthogonal complement of $E_{\sigma k}$; this 
proves the claim.

Hence $ad(\xi)$ takes the same eigenvalues $\lambda_j-\lambda_k$ 
on $\so_n^c$ as on $\u_n^c$. Thus by C1, these differences are 
integers and by C2 we even have $\lambda_{j+1}-\lambda_j = 
1$ as before; otherwise $\so_n^c$ had to be contained in a subalgebra 
$\Hom(E,E) + \Hom(F,F) \subset \u_n^c$ for some nontrivial 
decomposition $\C^n = E \oplus F$, but the inclusion
$SO_n \subset U_n$ is an irreducible representation. Thus we conclude 
that the set of eigenvalues $\lambda_j$ of ${1\over i}\xi$ is of 
the form $\{-r,-r+1,...,r-1,r\}$ for some positive integer or 
half integer $r$. Relabelling $E_j$ we obtain:

\theorem{Proposition A2.} An element $\xi \in \g = \so_n$ is 
canonical if and only if $\xi = \sum_{j=-r}^r j\.E_j$ for some 
orthogonal decomposition $\C^n = \sum_{j=-r}^r E_j$ such 
that $E_{-j} = \overline{E_j}$ for all $j \in \{-r,...,r\}$, for 
some $r \in {1\over2}\,\N$. Then $\g_k = \sum_j \tilde H_{j,j+k}$ 
where $\tilde H_{j,l} = \{A \in \hat H_{j,l};\ A^T = -A\}$.
\endtheorem

The corresponding symmetric space $S$ is a subset of the complex 
Grassmannian obtained from $E_{ev}$, namely $S = \{A(E_{ev});\ A 
\in SO_n\}$, where $E_{ev} := \sum_{j+r\ {\rm even}} E_j$. We have 
to distinguish two cases: 

\smallskip
\goodbreak
\item{(a)}  $r \in \N$: Then the eigenvalues of ${1\over i}
            ad(\xi)$ are the integers $j \in \{-r,...,r\}$. If 
            $j+r$ is even, then so is $-j+r$. Hence $E_{ev}$ is 
            invariant under conjugation and thus the 
            complexification of a subspace of $\R^n$. Hence $S$ is 
            the real Grassmannian containing all subspaces of 
            $\R^n$ with the same dimension as $E_{ev}$.

\smallskip
\item{(b)}  $r \not\in \N$: Then all eigenvalues $j \in 
            \{-r,...,r\}$ are proper half integers. If $j+r$ is even,  
            $-j+r$ is odd, and hence $\overline{E_{ev}} = 
            \sum_{j+r\ {\rm odd}} E_j = (E_{ev})^\perp$. Thus the 
            dimension $n$ is even and $E_{ev}$ is a maximal 
            isotropic subspace. Therefore $S$ is the space of all 
            maximal isotropic subspaces of $\C^n$, or as a coset 
            space, $S = SO_n/U_{n/2}$.

\theorem{Corollary.} The flag manifolds over real Grassmannians 
are precisely the manifolds of all orthogonal decompositions 
$\C^n = \sum_{j=-r}^r E_j$ for some $r \in \N$, where $E_{-j} = 
\overline{E_j}$ and the dimensions of $E_0,...,E_r$ are fixed 
arbitrarily.
\endproof
\endtheorem

\bigskip
The complexified tangent space of a general canonically embedded
flag manifold $Z = Ad(G)\xi$ at the point $\xi$ is
$
          T^c = ad(\g^c)\xi = ad(\xi)(\sum_j \g_j) 
                     = \sum_{j\neq 0} \g_j.  
$
Moreover, $Z$ is also a complex manifold (a coset space of the 
complex group $G^c$), and the space of (1,0) tangent vectors is $T' = 
\sum\limits_{j>0} \g_j$. Further, the complexified horizontal subspace 
for the fibration $\pi : Z \to S$ is ${\cal H} = \sum\limits_{k\ odd} 
\g_k$ while the (1,0) superhorizontal space is just ${\cal H}_1' 
= \g_1 \subset \cal H$.

In particular, for a flag manifold $Z$ over a real Grassmannian 
we obtain using the previous notation:
$$
             T^c = \sum_{j\neq k} \tilde H_{jk}, \ \ 
               T' = \sum_{j<k} \tilde H_{jk}, \ \
              {\cal H}_1' = \sum_j \tilde H_{j,j+1}.
\leqno\ (A5)
$$

\bigskip
\bigskip
\goodbreak
\leftline{\mittel References.}
\bigskip

\lit{[APS]} C.Arezzo, G.P.Pirola, M.Solci: {\it The Weierstrass
            representation for pluriminimal submanifolds},
            Preprint Pavia 2001
\lit{[BR]}  F.E.Burstall, J.H.Rawnsley: {\it Twistor Theory for 
            Riemannian Symmetric Spaces}, Springer L.N. in Math. 
            {\bf 1424}, 1990
\lit{[DG]}  M.Dajczer, D.Gromoll: {\it Real K\"ahler submanifolds 
            and uniqueness of the Gauss map}, J. Diff. Geom. {\bf 
            22} (1985), 13 - 28
\lit{[EH]}  J.-H. Eschenburg, E. Heintze: {\it Extrinsic symmetric 
            spaces and orbits of s\--representations}, manuscripta 
            math. {\bf 88} (1995), 517 - 524; {\it Erratum:}
            manuscr. math. {\bf 92} (1997), 408
\lit{[ET1]} J.-H.Eschenburg, R.Tribuzy: {\it (1,1)-geodesic maps 
            into Grassmann manifolds}, Math. Z. {\bf 220} (1995),
            337 - 346
\lit{[ET2]} J.-H.Eschenburg, R.Tribuzy: {\it Associated families 
            of pluriharmonic maps and isotropy}, manuscripta 
            math. {\bf 95} (1998), 295 - 310
\lit{[FT]}  M.Ferreira, R.Tribuzy, {\it K\"ahlerian submanifolds 
            of $\R^n$ with pluriharmonic Gauss map}, Bull. Soc. Math. 
            Belg. {\bf 45} (1993), 183 - 197
\lit{[F]}   D.Ferus, {\it Symmetric submanifolds of Euclidean 
            space}, Math. Ann. {\bf 247} (1980), 81 - 93
\lit{[RT]}  M.Rigoli, R.Tribuzy: {\it The Gauss map for 
            K\"ahlerian submanifolds of $\R^N$}, Trans. Amer. 
            Math. Soc. {\bf 332} (1992), 515 - 528
\lit{[RV]}  E.Ruh, J.Vilms: {\it The tension field of the Gauss 
            map}, Trans. Amer. Math. Soc. {\bf 149} (1970), 569 - 
            573
\lit{[Sp]}  M.Spivak: {\it A Comprehensive Introduction to 
            Differential Geometry, vol. IV}, Publish or Perish 
            1975
\lit{[Y]}   S.T.Yau: {\it Submanifolds with constant mean 
            curvature, I}, Amer. J. of Math. {\bf 96} (1974), 346 
            - 366

\bigskip
\bigskip
{\tt
f.e.burstall@maths.bath.ac.uk

eschenburg@math.uni-augsburg.de

mjferr@lmc.fc.ul.pt

tribuzy@buriti.com.br
}

\bye